\documentclass[twoside]{article}    

\usepackage[all]{xy}

\title{Varieties of Simultaneous Sums of Powers for binary forms}  

\author{Enrico Carlini \footnote{the author is a member of the PRIN01 project ``Spazi di moduli e teoria di Lie'' of the University of Pavia and a fellow of the OMATS programme of the University of Oslo.}}         

\date{}

\begin{document}           

\maketitle                 

\hyphenation{ ab-bia-mo al-ge-bri-ca-men-te as-so-cia-to as-so-cia-ti ca-rat-te-riz-za-zio-ne cer-chia-mo con-clu-de-re con-fi-gu-ra-zio-ni con-si-de-ra-re con-si-de-ria-mo co-ef-fi-cien-ti co-strui-re de-di-car-ci de-fi-nia-mo de-ri-van-ti di-cia-mo dif-fe-ren-zia-bi-le di-men-sio-ne di-men-sio-na-le di-mo-stra-zio-ne di-mo-stra-zio-ni di-su-gua-glian-ze e-qui-va-len-ti ge-ne-ra-to-ri e-le-men-ta-ri e-stre-ma-le fi-ni-ta-men-te fun-zio-ne ge-ne-ra-no ge-ne-ra-to ge-ne-ri-co gra-dua-to gra-dua-zio-ni li-near-men-te i-mi-ta-re im-me-dia-to in-cro-cio in-di-chia-mo in-di-spen-sa-bi-li in-di-stin-gui-bi-li i-po-te-si mo-stria-mo o-pe-ra-zio-ne o-pe-ra-zio-ni os-ser-via-mo ot-te-nia-mo pre-ce-den-te pa-ral-le-li-smi per-met-te pos-sia-mo pro-se-gui-re quan-do rap-pre-sen-ta-no re-sti-tu-isce ri-cor-da-re ri-go-ro-so si-gni-fi-ca-ti-vi sol-le-va-bi-le sor-go-no stu-dia-re suc-ces-sio-ni suf-fi-cien-te sup-po-nia-mo tra-dur-re va-rie-ta Za-ri-ski}


   \font\elevenmsb=msbm10 at 11pt            
   \font\eightmsb=msbm8
   \font\twelvemsb=msbm10 at 14pt 
   \font\twentysb=msbm10 at 21pt
   \newfam\msbfam
   \textfont\msbfam=\elevenmsb
   \scriptfont\msbfam=\eightmsb
   \scriptscriptfont\msbfam=\eightmsb
   \newcommand{\Bbb}[1]{{\fam\msbfam #1}}


  \newcommand{\GG}{\mbox{$\Bbb G$}}


  \newcommand{\PP}{\mbox{$\Bbb P$}}


  \newcommand{\Nat}{\mbox{$\Bbb N$}}


  \newcommand{\Integ}{\mbox{$\Bbb Z$}}


  \newcommand{\Cmp}{\mbox{$\Bbb C$}}


  \newcommand{\dimC}{\mbox{$\dim_{\Bbb C}$}}


  \newcommand{\rk}{\mbox{rk }}


  \newcommand{\LinSub}[1]{\mbox{LinSub}\left( #1 \right)}


  \newcommand{\LinSubI}[1]{\mbox{LinSub}_{#1}\left( H \right)}


  \newcommand{\SubD}[1]{\mbox{Sub}_{#1}\left( H \right)}



  \newcommand{\PolyR}[1]{k[X_1,\ldots,X_{#1}]}

  \newcommand{\PPolyR}[1]{k[X_0,\ldots,X_{#1}]}

  \newcommand{\bfv}{{\underline{v}}}


  \newcommand{\al}{\alpha}
  \newcommand{\sig}{\sigma}


  \newcommand{\calA}{{\cal A}}
  \newcommand{\calB}{{\cal B}}
  \newcommand{\calC}{{\cal C}}
  \newcommand{\calV}{{\cal V}}

  \def\TV{{\bf TV}}

  \def\HF{{\bf HF}}

  \def\TypeT{\mbox{$\cal T$}}
  \def\typeT{{\cal T}} 

  \newcommand{\TdotT}[1]{(\typeT_1,\ldots,\typeT_{#1})}

  \newcommand{\Ntype}[1]{\mbox{$ #1 $-}tipo}

  \newcommand{\Ntypes}[1]{\mbox{$ #1 $-}tipi}

  \newcommand{\OneType}[1]{\Bigl(\:(#1)\:\Bigr)}

  \newcommand{\TwoType}[2]{\Bigl(\:(#1) ,\,(#2)\:\Bigr)}

  \newcommand{\ThreeType}[3]{\Bigl(\:(#1) ,\,(#2),\,(#3)\:\Bigr)}

  \newcommand{\FourType}[4]{\Bigl(\:(#1) ,\,(#2),\,(#3),\,(#4)\:\Bigr)}



  \newcommand{\NTtype}[2]{\mbox{$(#1,#2)$}-tipo}
  \newcommand{\NTtypes}[2]{\mbox{$(#1,#2)$}-tipi}


  \newcommand{\NTtypeT}[2]{\typeT^{#1}_{#2}}


  \newcommand{\NTVtypeT}[3]{\typeT^{#1,{#3}}_{#2}}
 

  \newcommand{\VtypeT}[1]{\typeT^{#1}}

  \newcommand{\NTVh}[3]{h^{#1,{#3}}_{#2}}

  \newcommand{\Proj}[1]{{\Bbb P}^{#1}}


  \newcommand{\kProj}[1]{{\Bbb P}^{#1}_k}


 \newcommand{\kAff}[1]{{\Bbb A}^{#1}_k}
 \newcommand{\Aff}[1]{{\Bbb A}^{#1}}


  \newcommand{\pX}{{\Bbb X}}
  \newcommand{\ptY}{{\Bbb Y}}
  \newcommand{\ptZ}{{\Bbb Z}}


  \newcommand{\ptX}{{\Bbb X}}


  \newcommand{\Lin}[1]{{\Bbb #1}}
  \newcommand{\indLin}[2]{{\Bbb #1}_{#2}}


  \newcommand{\HilbS}[1]{\mbox{${\cal S}_{#1}$}}



  \def\sqw{\hbox{\rlap{\leavevmode\raise.3ex\hbox{$\sqcap$}}$%
\sqcup$}}
\def\sqb{\hbox{\hskip5pt\vrule width4pt height6pt depth1.5pt%
\hskip1pt}}
\def\cqfd{\ifmmode\sqw\else{\ifhmode\unskip\fi\nobreak\hfil
\penalty50\hskip1em\null\nobreak\hfil\sqw
\parfillskip=0pt\finalhyphendemerits=0\endgraf}\fi}


  \newcommand{\corrI}{\mbox{\boldmath $I$}}
  \newcommand{\corrV}{\mbox{\boldmath $V$}}

  \newcommand{\Hilb}{Hilbert}


  \newcommand{\card}[1]{\mid\!{#1}\!\mid}


  \newenvironment{Def}[1]{\par\medskip\noindent{\bf Definition\hskip 0.4cm }#1}{\par\medskip}

  
  \newtheorem{Teo}{Theorem}[section]

  
  \newtheorem{TeoDef}[Teo]{Teorema-Definizione}


  \newtheorem{Prop}[Teo]{Proposition}


  \newtheorem{PropDef}[Teo]{Proposition - Definizione}


  \newtheorem{Lem}[Teo]{Lemma}


  \newtheorem{Cor}[Teo]{Corollary}


  \font\Dimfont=cmr8            
  \newenvironment{Proof}[1]{\par\smallskip\noindent{P{\Dimfont ROOF}\hskip .5cm}#1}{\begin{flushright}\cqfd\end{flushright}\smallskip}


  \font\Dimfont=cmr8            
  \newenvironment{Dim}[1]{\par\smallskip\noindent{D{\Dimfont IMOSTRAZIONE}\hskip .5cm}#1}{\begin{flushright}\qed\end{flushright}\smallskip}

  \newtheorem{ProEs}[Teo]{Esempio} 
  \newenvironment{Es}[1]{\begin{ProEs}{\em #1}\vskip .5 cm \end{ProEs}}{}

  \newtheorem{ProRemark}[Teo]{Remark} 
  \newenvironment{Remark}[1]{\begin{ProRemark}{\em #1}\vskip .5cm \end{ProRemark}}{}

\begin{abstract}
The problem of simultaneous decomposition of binary forms as sums of powers of linear forms is studied. For generic forms the minimal number of linear forms needed is found and the space parametrizing all the possible decompositions is described. These results are applied to the study of rational curves.
\end{abstract}

\section{Introduction}  

Consider the polynomial ring $S=K[x_0,\ldots,x_n]$ and a form $f\in S_d$. A well known problem deals with the possible decompositions of  $f$ as a sum of powers of linear forms, that is
\[f=c_1l_1^d+\ldots +c_kl_k^d\]
$l_i\in S_1$, $c_i\in K$.

In geometric terms the problem reads as follow: given a point $[f]\in\PP S_d$ find points $[l_1^d],\ldots,[l_k^d]$ on the Veronese variety $\nu_d(\PP S_1)$ such that the $k$-secant space they span contains $[f]$.

We can ask for the minimal number of linear forms needed to decompose a given $f\in S_d$. Define this number as
\[k_{\min}(f)=\min\left\lbrace\,k: \exists\exists l_1,\ldots,l_k\in S_1,\quad f=\sum^k_{i=1}c_i\; l_i^d\,\right\rbrace .\]
When a decomposition exists this is usually not unique, so it is interesting to study all possible ways of decomposing a form $f\in S_d$ using exactly $k$ linear forms. With this in mind we set
\[VSP(f;k)=\overline{\left\lbrace\,\lbrace L_1,\ldots,L_k\rbrace\in\mbox{Hilb}_k\check{\PP S_1}: L_i\neq L_j, \quad f=\sum^k_{i=1}c_i l_i^d\,\right\rbrace}\]
where $L_i$ is the hyperplane of $\PP S_1$ defined by $l_i=0$ and $\mbox{Hilb}_k\check{\PP S_1}$ is the Hilbert scheme of length $k$ subschemes of $\check{\PP S_1}$.

Both $k_{\min}(f)$ and $VSP(f,k)$ were classically studied, but much remains unknown about them. An expected value for $k_{\min}$ was obtained by a naive parameter count, but  only recently it was proved that this value is exact for a generic form $f$, see \cite{AH}. The study of the variety $VSP(f,k)$ is still challenging and only few results are known in  general. For more on this see \cite{RS}.

A straightforward generalization is the study of the simultaneous decompositions of a set of forms $f_1,\ldots,f_r\in S_d$, that is
\[f_i=c_{i1}l_1^d+\ldots c_{ik}l_k^d,\quad i=1,\dots,r\]
involving the same linear forms $l_j$.

In this case also we have a geometric interpretation: given points $[f_1],\dots,[f_r]\in\PP S_d$ find points $[l_1^d],\ldots,[l_k^d]$ on the Veronese variety $\nu_d(\PP S_1)$ such that the $k$-secant space they span contains the linear space $<[f_1],\dots,[f_r]>$.

This problem was classically studied by means of polar polyhedra, e.g. see \cite{L} for $n=2$ and $d=3$,  and in a more general setting  by Terracini. In \cite{T} a solution for the case $n=2$, $r=2$ is claimed and a general criterion is stated. For a  rigorous proof and generalization of Terracini's result see \cite{DF}. For an exposition in modern terms and an interesting interpretation of Terracini's criterion  see \cite{CC}.

As in the case of one form, there are two main objects of interest.
\begin{Def} Let $f_1,\ldots,f_r\in S_d$. We define
\[k_{\min}(f_1,\ldots,f_r)=\min\left\lbrace\,k: \exists\exists l_1,\ldots,l_k\in S_1,\quad f_i=\sum^k_{j=1}c_{ij}\; l_j^d\quad i=1,\dots,r\,\right\rbrace .\]
\end{Def}

\begin{Def}
Let $f_1,\ldots,f_r\in S_d$. We define the {\bf v}ariety of {\bf s}imultaneous {\bf s}ums of {\bf p}owers of the $f_i$'s with respect to $k$ to be
\[VSSP(f_1,\ldots,f_r;k)=\overline{\left\lbrace\, \lbrace\,L_1,\dots,L_k\,\rbrace\in\mbox{Hilb}_k\check{\PP S_1}: L_i\neq L_j ,f_i=\sum_{j=1}^k c_{ij}\;l_j^d \quad i=1,\dots,r\,\right\rbrace} ,\]
where $L_i$ is the hyperplane of $\PP S_1$ defined by $l_i=0$ and $\mbox{Hilb}_k\check{\PP S_1}$ is the Hilbert scheme of length $k$ subschemes of $\check{\PP S_1}$.
\end{Def}
We notice that if  $f_1,\dots,f_r$ are linearly dependent, then the problem reduces to the study of $r'$ independent forms, $r'<r$. Therefore we may assume the $f_i$'s to be linearly independent.

As in the $r=1$ case there is an expected value for $k_{\min}$ obtained by a parameter count:
\[\min\left\lbrace\,{r \over r+n}{n+d \choose d}, \dim G(r,S_d)\,\right\rbrace ,\]
where $G(r,S_d)$ is the Grassmannian of $r$ dimensional subspaces of $S_d$.\hfill\break
When $r>1$ only few values of $(n,d,r,k)$ are known for which $k_{\min}$ is not the  predicted value for a generic choice of forms, see \cite{CC}. These exceptions can be proved by ad hoc methods but there are few general results. The most general result about $k_{\min}$  asserts that, in the binary case $n=1$, the actual and the expected value of $k_{\min}(f_1,\dots,f_r)$  are equal  for generic forms. This is proved in \cite{CC} but the proof is quite involved. To our knowledge there are  almost no results  about the variety $VSSP$ in the case $r>1$.

In this paper we  restrict our attention to the binary case $S=K[x_0,x_1]$. In this situation we find a formula for $k_{\min}(f_1,\ldots,f_r)$ in a more direct way than in \cite{CC}. Moreover we give a complete description of $VSSP(f_1,\ldots,f_r;k)$ for $f_1,\ldots,f_r$ generic forms. In the last section we  show an application of $VSSP$ to the study of rational curves.

This problem was proposed to the author by Prof. Kristian Ranestad during the international school PRAGMATIC 2000 held in Catania. The author wishes to thank the OMATS programme for the financial support during the preparation of this work. The author also wishes to give a special acknowledgment to Prof. Kristian Ranestad for his many fruitful suggestions and comments.

\section{Apolarity and inverse systems}
Set $S=K[x_0,\ldots,x_n]$ and $T=K[y_0,\ldots,y_n]$ with $K=\overline K$ a field of characteristic 0. We make $S$ into a $T$-module via differentiation: given monomials $y^\alpha$, $x^\beta$ we define
\[y^\alpha \circ x^\beta =\left\lbrace\begin{array}{lr}0 & \mbox{if $\alpha_i>\beta_i$ for some $i$} \\ \alpha!{\beta  \choose \alpha }x^{\beta-\alpha}  & \mbox{otherwise}\end{array}\right. ,\]
where the computations on the multi-indices are made componentwise, e.g. $\alpha!=\alpha_0!\ldots\alpha_n!$.

We recall that:
\begin{itemize}
\item given $f\in S$ we set $f^\perp=\lbrace\,D\in T : D\circ f=0\,\rbrace$, this is an ideal in $T$ and it is called the {\em orthogonal ideal of \/}$f$;
\item given $D\in T$ we set $D^{-1}=\lbrace\,f\in S : D\circ f=0\,\rbrace$, this is a graded $T$-submodule of $S$ and it is  called the {\em inverse system of\/}$D$.
\end{itemize}
We will need some basic properties about orthogonal ideals and of inverse systems:

\noindent{\bf Properties} (see \cite{G}, pp 11-19)
\begin{enumerate}
\item\label{Gorperp} if  $f\in S_d$,  then $f^\perp$ is a Gorenstein artinian ideal with socle degree $d$;
\item\label{DinverseDIM} if $D\in T_k$, $d\geq k$, then $\dim_K (D^{-1})_d={n+d \choose n}-d+k-1$;
\item\label{D-1inj} if $D,G\in T_k$, $d\geq k$, such that $(D^{-1})_d= (G^{-1})_d$, then $[D]=[G]$ in $\PP T_k$;
\item\label{DinversePERP} if $f_1,\ldots,f_r\in S_d$, then
\[\left(\bigcap_{i=1}^r f_i^\perp\right)_k=\lbrace\,D\in T_k : (D^{-1})_d\supset <f_1,\ldots,f_r>\,\rbrace ;\]
\item\label{dual} via apolarity we have a natural identification $\check{\PP S_k}=\PP T_k$.
\end{enumerate}

Apolarity is a powerful tool in studying the decomposition of forms as sums of powers because of the following (see \cite{RS}, 1.3)
\begin{Lem}[Apolarity Lemma]\label{apolaritylemma}
Let $f_1,\ldots,f_r\in S_d$,  then the following facts are equivalent:
\begin{enumerate}
\item $\exists\exists\,  c_{ij}\in K$, $l_1,\ldots,l_k\in S_1$, $[l_a]\neq [l_b]$ in $\PP S_1$ for $a \neq b$, such that \[f_i=\sum_{j=1}^k c_{ij}\;l_j^d, \quad i=1,\dots,r ;\]
\item $\exists\exists\,  L_1,\ldots,L_k\in\PP T_1, L_i\neq L_j$ for $i\neq j$, such that
\[I_\Gamma\subset\bigcap_{i=1}^r f_i^\perp\]
 where $I_\Gamma$ is the ideal of the set of points $\Gamma=\lbrace\,L_1,\ldots,L_k\,\rbrace$.
\end{enumerate}
\end{Lem}

The lemma leads our attention to the ideals contained in $\bigcap_{i=1}^r f_i^\perp$, so we give the following
\begin{Def}
Given $f_1,\ldots,f_r\in S_d$ the {\bf v}ariety of {\bf s}imultaneous a{\bf p}olar {\bf s}usbschemes of length $k$ of the $f_i$'s is
\[VSPS(f_1,\ldots,f_r;k)=\left\lbrace\,\Gamma\in\mbox{Hilb}_k\check{\PP S_1} : I_\Gamma\subset \bigcap_{i=1}^r f_i^\perp\,\right\rbrace .\]
\end{Def}

The Apolarity Lemma shows why the binary case is easier to treat. When $n=1$ the ideal of  a set of points is a principal ideal and the generator is square free  if the points are distinct. Hence there is a natural identification
\[\mbox{Hilb}_k\check{\PP S_1}=\check{\PP S_k}=\PP T_k,\]
where the last equality comes from Property \ref{dual}.

Finally, using Property \ref{DinversePERP} and Lemma \ref{apolaritylemma}, we get a useful description of $VSPS$ and of $VSSP$. Let $S=K[x_0,x_1]$ and $f_1,\ldots,f_r\in S_d$. Then\[\]
\hfill$\begin{array}{lcl}VSPS(f_1,\ldots,f_r;k) & = & \lbrace\,[D]\in\PP T_k : D\in(\cap_{i=1}^rf_i^\perp)_k \,\rbrace \\ \\ & = & \PP\left(\bigcap_{i=1}^rf_i^\perp\right)_k \\ \\ &=& \lbrace\,[D]\in\PP T_k : (D^{-1})_d\supset <f_1,\ldots,f_r>\,\rbrace  \end{array}$\hfill$(\dagger)$\break
\[\]
and
\[\]
\hfill$\begin{array}{lcl}VSSP(f_1,\ldots,f_r;k) & = & \overline{\lbrace\,[D]\in\PP T_k : D\in(\cap_{i=1}^rf_i^\perp)_k, D\not\in\Delta_k \,\rbrace} \\ \\ & = &  \PP\overline{\left(\bigcap_{i=1}^rf_i^\perp\right)_k\setminus\Delta_k }\\ \\ &=& \overline{\lbrace\,[D]\in\PP T_k : (D^{-1})_d\supset <f_1,\ldots,f_r>,D\not\in\Delta_k\,\rbrace} \end{array}$\hfill$(\ddagger)$\[\]
\noindent where $\Delta_k$ is the discriminant locus of polynomials of degree $k$ with at least a repeated root. We also notice that
\[k_{\min}(f_1,\ldots,f_r)=\min\lbrace\,k : \exists D\in T_k : (D^{-1})_d\supset <f_1,\ldots,f_r>,D\not\in\Delta_k \,\rbrace .\]

It is useful to summarize these results.
\begin{Prop}\label{bigcap}
If $S=K[x_0,x_1]$ and $f_1,\ldots,f_r\in S_d$, then  $VSPS(f_1,\ldots,f_r;k)$ and $VSSP(f_1,\ldots,f_r;k)$ are projective spaces for every $k$. Moreover 
\[VSPS(f_1,\ldots,f_r;k)\supseteq VSSP(f_1,\ldots,f_r;k)\]
and they are  equal whenever the latter is not empty.
\end{Prop}

Given explicit binary forms $f_1,\ldots,f_r$ we can actually determine $k_{\min}$, $VSSP$ and $VSPS$ using $(\dagger)$ and $(\ddagger)$. This requires an easy algorithm involving linear algebra (orthogonal ideals) and basic Gr\"obner basis computations (intersection of ideals).\hfill\break
The really tough problem is deriving results holding for a generic choice of $r$ forms. Part of the difficulty is related to the bad behavior of orthogonal ideals. For example it easy to show that for any binary form $f\in S_d$
\[\dim_K (f^\perp)_k\geq 2k-d,\]
but the actual value of the dimension depends on the particular form we choose.

The best result we can obtain for orthogonal ideals is an easy consequence of the previous bound and of Grassmann's formula for vector spaces.

\begin{Lem}\label{grassmann}
Let $d,r,k$ be natural numbers and $S=K[x_0,x_1]$. If

\[k>{r(d+1)-1\over r+1},\]

then for any choice of $f_1,\ldots,f_r\in S_d$ we have

\[(f_1^\perp\cap\ldots\cap f_r^\perp)_k\neq (0).\]

\end{Lem}

\section{The geometric setting}
From now on  we will consider only binary forms, so that $S=K[x_0,x_1]$ and $T=K[y_0,y_1]$.

Consider the map
\[\begin{array}{llcl}
\psi_k: & \PP T_k & \longrightarrow & G(k, S_d) \\
        & [D]     & \mapsto & (D^{-1})_d
\end{array};\]
by Property \ref{DinverseDIM} it is well defined when $k\leq d$. Using Pl\"{u}cker coordinates and Property \ref{D-1inj} one  verifies that $\psi_k$ is an isomorphism on its image, $G_k=\psi(\PP T_k)$, for all $k\leq d$. If we let $\Delta_k\subset T_k$ be the locus of forms with repeated roots, then we have
\[ \dim G_k=\dim G_{\Delta_k}+1,\]
where $G_{\Delta_k}=\psi_k(\PP\Delta_k)$.

Now consider the following diagram


\[\begin{array}{rcll}
G(r, S_d)\times G(k, S_d)  & \supset    & \Sigma_k=\lbrace\,(\Lambda,\Gamma):\Lambda\subseteq\Gamma,\Gamma\in G_k\,\rbrace & \supset\Sigma_{\Delta_k}=\lbrace\,(\Lambda,\Gamma):\Lambda\subseteq\Gamma,\Gamma\in G_{\Delta_k}\,\rbrace \\
                          &  &\downarrow\varphi_k & \\
G(r,S_d)          & \supset    & \widetilde\Sigma_k & 
\end{array}\]

The study of the simultaneous decompositions of a set of forms $f_1,\ldots,f_r$ is equivalent to the study of the map $\varphi_k$, as shown by the following
\begin{Prop}\label{geoversionofprb}
Let $f_1,\ldots,f_r\in  S_d$ be linearly independent forms and let  $\Lambda=<f_1,\ldots,f_r>$. Then
\begin{enumerate}
\item $k_{min}(f_1,\ldots,f_r)=\min\lbrace\,k: \varphi_k^{-1}(\Lambda)\setminus \Sigma_{\Delta_k}\neq\emptyset\,\rbrace$;
\item $VSSP(f_1,\ldots,f_r;k)\simeq\overline{\varphi_k^{-1}(\Lambda)\setminus\Sigma_{\Delta_k}}$;
\item $VSPS(f_1,\ldots,f_r;k)\simeq\varphi_k^{-1}(\Lambda)$.
\end{enumerate}
\end{Prop}
\begin{Proof}{
First we compute the fiber of $\varphi_k$ on $\Lambda$:
\[\varphi_k^{-1}(\Lambda)=\lbrace\,(\Lambda,(D^{-1})_d):\Lambda\subset(D^{-1})_d, [D]\in\PP T_k\,\rbrace;\]
from this we immediately get part {\it 1.}

Now, using $(\ddagger)$,  we notice that
\[\psi_k^{-1}(\overline{\varphi_k^{-1}(\Lambda)\setminus\Sigma_{\Delta_k}})=\overline{\lbrace\,[D]\in\PP T_k:\Lambda\subset(D^{-1})_d, D\not\in\Delta_k \,\rbrace}=\PP\overline{(\cap_if_i^\perp)_k\setminus\Delta_k}.\]
Because $\psi_k$ is an isomorphism we get part  {\it 2}. The same argument and $(\dagger)$ give part {\it 3.}
}\end{Proof}

The map $\varphi_k$ is useful in solving our problem also because of the properties of the varieties $\Sigma_k$ and $\Sigma_{\Delta_k}$. In fact we have

\begin{Lem}\label{sigmairr}
$\Sigma_k$ and $\Sigma_{\Delta_k}$ are Grassmannian bundles on irreducible varieties. In particular they are irreducible  and we have
\[\dim \Sigma_k= \dim\Sigma_{\Delta_k}+1=k+r(k-r). \]
\end{Lem}
\begin{Proof}{
It is enough to consider the projection maps
\[\Sigma_k\longrightarrow G_k\]
\[\Sigma_{\Delta_k}\longrightarrow G_{\Delta_k}\]
and to notice that their fibers are the Grassmannians $G(r,k)$.
}\end{Proof}

Finally we can state our main result.
\begin{Teo}\label{main}Let $S=K[x_0,x_1]$. Given natural numbers $d,r$  set \[k_{\min}(d,r)=\min\left\lbrace\,k:k>{r(d+1)-1\over r+1}\,\right\rbrace .\]
There exists an open non-empty subset
\[V_{d,r}\subset G(r,S_d)\]
such that, for all $f_1,\ldots,f_r\in S_d$ satisfying $[<f_1,\ldots,f_r>]\in V_{d,r}$, the following hold:
\begin{enumerate}
\item $k_{\min}(f_1,\ldots,f_r)=k_{\min}(d,r)$;
\item $VSSP(f_1,\ldots,f_r;k)=\left\lbrace\begin{array}{lc}\PP^{k(r+1)-r(d+1)} & k\geq k_{\min}(d,r)\\ \\ \emptyset & k<k_{\min}(d,r)\end{array}\right.$;
\item $VSPS(f_1,\ldots,f_r;k)=\emptyset$ if $ k<k_{\min}(d,r)$.
\end{enumerate}
Moreover
\[VSSP(f_1,\ldots,f_r;k_{\min}(d,r))=\left\lbrace
\begin{array}{ll}
\PP^{r+1-\varepsilon} & \varepsilon\neq 0 \\ \\
\lbrace\mbox{pt}\rbrace & \varepsilon = 0 
\end{array}\right.\]
where $\varepsilon\equiv r(d+1)$ {\em mod}$(r+1)$.
\end{Teo}
\begin{Proof}{
Set $Z_k=\lbrace\,\Lambda\in\widetilde\Sigma_k : \varphi_k^{-1}(\Lambda)\subset\Sigma_{\Delta_k}\,\rbrace$ and consider the following diagram
\[\begin{array}{lcccc}
\Sigma_k & \supset & \Sigma_{\Delta_k} & \supset & \varphi_k^{-1}(Z_k) \\
 \downarrow\varphi_k& & & \\
\widetilde\Sigma_k & \supset & Z_k &
\end{array}\]
Set $\lambda=\min\lbrace\,\dim\varphi_k^{-1}(\Lambda):\Lambda\in\widetilde\Sigma_k\,\rbrace$. Using Lemma \ref{sigmairr} and the Fiber Dimension Theorem (see \cite{H}, lecture 11) we get
\[\dim\Sigma_k=\dim\widetilde\Sigma_k+\lambda,\]
\[\dim\Sigma_{\Delta_k}\geq\dim Z_k+\lambda.\]
Hence $\dim\widetilde\Sigma_k-\dim Z_k\geq 1$.

If $\widetilde\Sigma_k$ is dense, then $\dim\Sigma_k\geq\dim G(r,S_d)$. This gives the condition
\[k\geq k_1=\min\left\lbrace\, k : k\geq {r(d+1)\over r+1}\,\right\rbrace .\]
By Lemma \ref{grassmann} we know that if $k\geq k_2=\min\lbrace\,k:k>{r(d+1)-1\over r+1}\,\rbrace$ then $\varphi_k$ is dense.
It is easy to check that $k_1=k_2=k_{\min}(d,r)$. For the sake of simplicity set  $\bar k=k_{\min}(d,r)$.

Finally we set
\[U=\widetilde\Sigma_{\bar k}\setminus(\overline{Z_{\bar k}\cup\widetilde\Sigma_{\bar k-1}}).\]
By the preceding consideration $U\subset G(r, S_d)$ is open and non-empty. Moreover, if $\Lambda\in U$ then
\[\varphi_{\bar k}^{-1}(\Lambda)\setminus\Sigma_{\Delta_{\bar k}}\neq\emptyset,\]
\[\varphi_{ k}^{-1}(\Lambda)=\emptyset, k<\bar k.\]
Using Proposition \ref{geoversionofprb} we conclude that $U$ satisfies part {\it 1.}

Proposition \ref{bigcap}  yields
\[VSSP(f_1,\ldots,f_r;k)=VSPS(f_1,\ldots,f_r;k)=\emptyset,\quad k<\bar k=k_{\min}(d,r)\]
for $f_1,\ldots,f_r\in S_d$ such that  $<f_1,\ldots,f_r>\in U$. This proves part {\it 3}.

By Propositions \ref{bigcap} and \ref{geoversionofprb} we know that $VSSP(f_1,\ldots,f_r;k)$, $k\geq k_{\min}(d,r)$,  is a projective space of dimension $\dim\varphi_{ k}^{-1}(<f_1,\ldots,f_r>)$. As the fiber dimension is an upper  semicontinuous function,  there exists an open non-empty subset $U'\subset U$ such that
\[\dim VSSP(f_1,\ldots,f_r;k)=\dim\Sigma_k-\dim G(r,S_d)=k(r+1)-r(d+1)\]
for $<f_1,\ldots,f_r>\in U'$. This completes the proof of part {\it 2}.

To get the expression for $\dim VSSP(f_1,\ldots,f_r;k_{\min}(d,r))$ we only have to use part {\it 2} and to study the congruence class of $r(d+1)-1$ mod $r+1$.

Letting $V_{d,r}=U'$ completes the proof.
}\end{Proof}

{\bf Example} Using  Theorem  \ref{main} we can recover a classical result of Sylvester (1851). Given a generic binary form $f\in S_d$, i.e. $f\in V_{d,1}$, the minimal number of linear forms needed to decompose it is
\[k_{min}(f)= \min\lbrace\,k:k>{d \over 2}\,\rbrace=\left[{d \over 2}\right]+1\]
and the  possible  decompositions are parametrized by
\[VSSP(f,k_{min}(f))=\left\lbrace
\begin{array}{ll}
\PP^1 & d\mbox{ even} \\
\lbrace\mbox{pt}\rbrace & d\mbox{ odd} \\
\end{array}
\right..\]

\section{Rational curves}
\begin{Def}
A rational curve $C\in\PP^n$ is the image of a rational map $\alpha:\PP^1\rightarrow\PP^n$. We say that the curve is non-degenerate if it is not contained in a hyperplane.
\end{Def}

Let $S=k[x_0,x_1]$ and fix the standard lex ordered monomial basis, e.g. with respect to $x_0>x_1$, in each of the homogeneous pices $S_n$, so we have the identifications $\PP ^n=\PP S_n$ for all $n$.

Let $C\subset\PP^n$ be a rational curve of degree $d$, with $d>n$. There exists a unique $\Lambda_C\in\GG(d-n-1,\PP^d)=\lbrace\,\Lambda\subset\PP^d : \Lambda\simeq\PP^{d-n-1}\,\rbrace$ such that the following diagram commutes


\[\begin{array}{lccl}
 & & \PP^d & \supset C_d \\
 & \begin{array}{c}\nu_d \\ \nearrow \end{array}& \downarrow \pi & \\
\PP^1 & \longrightarrow & \PP^n & \supset C
\end{array}\]
where $\nu_d$ is the $d$-uple embedding of $\PP^1$, $C_d$ is the rational normal curve of degree $d$ and $\pi$ is the projection from $\Lambda_C$. In particular $\pi(C_d)=C$.
\begin{Def}
Let $C\subset\PP^n$ be a rational curve. We define
\[S_b^a(C)=\lbrace\,\Gamma\in\GG(a,\PP^n):\alpha^{-1}(\Gamma\cap C)\mbox{ has length }b\,\rbrace\]
\[S_b^a(C)_{\neq}=\lbrace\,\Gamma\in\GG(a,\PP^n):\alpha^{-1}(\Gamma\cap C)\mbox{ is smooth of length }b\,\rbrace.\]
\end{Def}
We notice that $S_b^0(C)$ is the set of $b$-uple points of $C$ and $S_b^0(C)_{\neq}$ is the set of $b$-uple points of $C$ having $b$ distinct tangent lines, e.g. if $C$ is a plane cubic with a node then $S_2^0(C)_{\neq}=\lbrace\mbox{pt}\rbrace$ while $S_2^0(C')_{\neq}=\emptyset$ if $C'$ is a cusp.

If $C$ is a smooth curve, then $S_b^a(C)_{\neq}$ has a  nicer geometric description:
\[S_b^a(C)_{\neq}=\lbrace\,\Gamma\in\GG(a,\PP^n):\Gamma\cap C \mbox{ is a set of $b$ distinct points}\,\rbrace.\]

\noindent Let $C\subset\PP^n$ be a rational non-degenerate curve of degree $d$.  It is immediate to verify the following:
\begin{itemize}
\item $S_d^{n-1}(C)=\check{\PP^n}$  as, taking multiplicities into account, {\it any} hyperplane intersects $C$ in $d$ points;
\item $S_d^{n-1}(C)_{\neq}$ is dense in $\check{\PP^n}$ as a {\it generic\/} hyperplane intersects $C$ in $d$ distinct points;
\item $S_{d^{'}}^{n-1}(C)=S_{d^{'}}^{n-1}(C)_{\neq}=\emptyset$ if $d^{'}\neq d$.
\end{itemize}

We notice that $S_b^a(C),S_b^a(C)_{\neq}$ are not interesting for all the values of $a$ and $b$, as shown by the following

\begin{Lem}\label{b-a}
Let $C\subset\PP^n$ be a rational non-degenerate curve of degree $d$, $d>n$. If $b-a> d-n+1$, then $S_b^a(C)=S_b^a(C)_{\neq}=\emptyset$.
\end{Lem}
\begin{Proof}{
Let $\Gamma\in S^a_b(C)$; then choosing $P\in C\setminus\Gamma$ we build $\Gamma^{'}=<\Gamma,P>\in S^{a+1}_{b^{'}}(C)$, $b^{'}\geq b+1$. Repeating the construction we get
\[\overline\Gamma\in S^{n-1}_{\overline b}(C)\]
where $\overline b\geq b+n-a-1$. As $S^{n-1}_{\overline b}(C)\neq\emptyset$  we must have $d-n+1\geq b-a$.
We have shown that $S^a_b(C)\neq\emptyset$ implies $d-n+1\geq b-a$ and by negation we get the thesis.
}\end{Proof}

We will describe $S_b^a(C),S_b^a(C)_{\neq}$ in the extremal case  $b-a= d-n+1$.
\begin{Prop}\label{curveVSSP}
Let $C\subset\PP^1$ be a rational  non-degenerate curve of degree $d$, $d>n$,  and $a,b$ natural numbers such that $b-a=d-n+1$. Then  the following isomorphisms hold
\[S_b^a(C)\simeq VSPS(f_1,\ldots,f_r;b)\]
\[S_b^a(C)_{\neq}\simeq V\subset VSSP(f_1,\ldots,f_r;b)\]
where $V$ is dense. 
\end{Prop}
\begin{Proof}{Let $\Lambda_C=<f_1,\ldots,f_r>$, $f_1,\ldots,f_r\in S_d$. As $C$ is the projection of $C_d$ from $\Lambda_C$ we have
\[\begin{array}{lcl}
S_b^a(C) & = & \lbrace\,\Gamma\in\GG(a,\PP^n):\alpha^{-1}(\Gamma\cap C)\mbox{ has length }b\,\rbrace \\
         & \simeq & \lbrace\,\Gamma\in\GG(a+d-n,\PP^n):\Gamma\supset\Lambda_C ,\nu_d^{-1}(\Gamma\cap C_d)\mbox{ has length }b\,\rbrace=A
\end{array}\]
and
\[\begin{array}{lcl}
S_b^a(C)_{\neq} & = & \lbrace\,\Gamma\in\GG(a,\PP^n):\alpha^{-1}(\Gamma\cap C)\mbox{ is smooth of length }b\,\rbrace \\
         & \simeq & \lbrace\,\Gamma\in\GG(a+d-n,\PP^n):\Gamma\supset\Lambda_C ,\nu_d^{-1}(\Gamma\cap C_d)\mbox{ is smooth of length }b\,\rbrace=B
\end{array}\]

First we construct the isomorphism $A\simeq VSPS(f_1,\ldots,f_r;b)$.\hfill\break
Given $\Gamma\in A$ consider the subscheme of $Y=\nu_d^{-1}(\Gamma\cap C_d)\subset \PP^1$: it is defined by the ideal
\[I_Y=(g_1,\ldots,g_{n-a})\]
where the $g_i$'s are the pullbacks of the linear equations defining $\Gamma$, as the scheme has length $b$ they have a common factor of degree $b$. The $g_i$'s are independent so that the numerical condition on $a$ and $b$ implies that the saturation of $I_Y$ is the principal ideal $(g)$, $g=GCD(g_i)\in S_b$. Let $G_i(y_0,y_1)=g_i(y_0,y_1)\in T_b$ and $G(y_0,y_1)=g(y_0,y_1)\in T_b$. By apolarity the hyperplanes $(G_i^{-1})_d$ define $\Gamma$. It is easy to verify that $(G^{-1})_d\supset\Gamma$, i.e. $G\in VSPS(f_1,\ldots,f_r;b)$, so we have a map $A\rightarrow VSPS(f_1,\ldots,f_r;b)$.\hfill\break
The inverse map is defined by sending $G\in VSPS(f_1,\ldots,f_r;b)$ to the linear space $(G^{-1})_d\in A$.

Now we construct the isomorphism $B\simeq V$, $V$ a dense subset of $VSSP(f_1,\ldots,f_r;b)$.\hfill\break
Let $V$ be the dense subset of $VSSP(f_1,\ldots,f_r;b)$ consisting of square free polynomials. Then we define the map $V\rightarrow B$ by sending $G\in V$ to the linear space $(G^{-1})_d$.\hfill\break
To define the inverse map, given $\Gamma\in B$,  we consider the subscheme $\nu_d^{-1}(\Gamma\cap C_d)$ which is defined by the principal ideal $(g),g\in S_b$: $G(y_0,y_1)=g(y_0,y_1)$ is a square free polynomial and $G\in VSSP(f_1,\ldots,f_r;b)$.
}\end{Proof}

{\bf Example} Let $f_1=-2x_0^5+2x_1^5+(x_0-x_1)^5$, $f_2=-6x_0^5+3x_1^5+2(x_0-x_1)^5$ and $\Lambda_C=<f_1,f_2>$. We want to study the rational curve $C=\pi(C_5)\in\PP^3$ obtained as projection of the rational normal curve of $\PP^5$ from $\Lambda_C$.\hfill\break
We know that
\[S^a_b(C)=S^a_b(C)_{\neq}=\emptyset\]
when $b-a>3$. Hence  the interesting cases are: 
\[S^a_{a+3}(C), S^a_{a+3}(C)_{\neq} \quad a=0,1,2 .\]
To apply Proposition \ref{curveVSSP}  we need to determine $I=f_1^\perp\cap f_2^\perp$. By direct computation we get
\[I=(y_0y_1(y_0+y_1),y_0^4+y_1^4)\cap(y_0y_1(y_0+y_1),y_0^4+2y_1^4),\]
in particular
\[\dim_K I_5=4,\quad \dim_K I_4=2, \quad \dim_K I_3=1, \quad \dim_K I_d=0\quad \mbox{for }d<3\]
It is possible to verify that each homogeneous piece of $I$ is generated by forms without common roots. Moreover the generator of $I_3$ is square free. Hence,  using $(\ddagger)$ and Proposition\ref{curveVSSP}, we obtain
\[S^0_3(C)_{\neq}=\PP I_3=\lbrace\mbox{pt}\rbrace,\]
\[S^1_4(C)_{\neq}=\PP I_4=\PP^1,\]
\[S^2_5(C)_{\neq}=\PP I_5=\PP^3.\]
We also have that $S^a_{a+3}(C)$ is dense in  $S^a_{a+3}(C)_{\neq}$ for $a=0,1,2$.\hfill\break
In particular \[S^0_3(C)_{\neq}=\lbrace\mbox{pt}\rbrace\] means that $C$ has a unique triple point.

We have shown how to use $VSSP$ to study curves that are projection of the rational normal curve: given the center of projection $\Lambda=<f_1,\ldots,f_r>$ we have to investigate the decompositions of the $f_1,\ldots,f_r$ as sums of powers of linear forms, so that each case has to be treated separately. To find general results we have to exclude curves with a pathological behavior and this can be done using Theorem \ref{main}.
\begin{Def}
A rational curve $C\in\PP^n$ of degree $d$, $d>n$, is said to be generic if $\Lambda_C\in V_{d,r}$, where $r=\dim_K \Lambda_C$ and $V_{d,r}$ is as given by Theorem \ref{main}.
\end{Def}
We notice that, because $V_{d,r}$ is open and non-empty, almost all the rational curves $C\subset\PP^n$ of degree $d$, $d>n$, are generic.

An easy consequence of Theorem \ref{main} and of Proposition \ref{curveVSSP} is the following:

\begin{Cor}\label{cor}
Let $C\subset\PP^n$ be a  generic rational curve of degree $d$, $d>n$. If  $b-a=d-n+1$ then, defining $k_{\min}(d,d-n)$ as in Theorem \ref{main}, we have
\begin{enumerate}
\item $S^a_b(C)=S^a_b(C)_{\neq}=\emptyset$ for $b<k_{\min}(d,d-n)$;
\item $S^a_b(C)=\PP^{b(d-n+1)-(d-n)(d+1)}$ for $b\geq k_{\min}(d,d-n)$;
\item $S^a_b(C)_{\neq}$ is dense in $\PP^{b(d-n+1)-(d-n)(d+1)}$ for $b\geq k_{\min}(d,d-n)$.
\end{enumerate}

\end{Cor}

{\bf Example} Let $C\subset\PP^3$ be a rational non-degenerate curve of degree 5. By Lemma \ref{b-a} we know that $S^a_b(C)=S^a_b(C)
_{\neq}=\emptyset$ for $b-a>3$.\hfill\break
If $C$ is generic, then using Corollary \ref{cor} and the fact that $k_{\min}(5,2)=4$  we get
\[S^0_3(C)_{\neq}=\emptyset\]
\[S^1_4(C)_{\neq}=\lbrace\mbox{pt}\rbrace\]
\[S^2_5(C)_{\neq}=\PP^3 .\]
In particular, because $S^0_3(C)=\emptyset$, we conclude that $C$  has no triple points. This shows that the curve of the previous example is not generic.\hfill\break
If $C$ is also smooth then the equality \[S^1_4(C)_{\neq}=\lbrace\mbox{pt}\rbrace\] means that there exists a unique 4-secant line to $C$.

{\bf Example} Let $C\subset\PP^{16}$ be a smooth generic rational  curve of degree 19. Because the curve is smooth we can get interesting geometric properties from studying $S^b_a(C)_{\neq}$.\hfill\break
$C$ is the projection of the rational normal curve $C_{19}$ from $\Lambda_C=<f_1,f_2,f_3>$ and, because the curve is generic,  we know that $k_{\min}(f_1,f_2,f_3)=15$. This means that
\[S^{a}_{a+4}(C)_{\neq}=VSSP(f_1,f_2,f_3; a+4)=\emptyset, \mbox{ for }a=0,\dots, 10.\]
We  get that \[S^{11}_{15}(C)_{\neq}=\lbrace\mbox{pt}\rbrace\] so that there exists a unique 15-secant $\PP^{11}$ to $C$. We also have 
\[S^{12}_{16}(C)_{\neq}=\PP^4,S^{13}_{17}(C)_{\neq}=\PP^8,S^{14}_{18}(C)_{\neq}=\PP^{12}.\]

\parbox{13cm}{\small 
Dipartimento di Matematica, Universit\`a degli Studi Pavia, Via Ferrata 1 27100 Pavia, Italy. \\ 
Matematisk Institutt, UiO, P.B. 1053 Blindern, N-0316 Oslo, Norway.\\
http:{\tt dimat.unipv.it/\mbox{ $\tilde{}$} carlini}\\
email: {\tt carlini@dimat.unipv.it}\\}

\end{document}